\documentstyle[11pt]{article}

\title{A SHORTER PROOF FOR THE TRANSITIVITY OF TRANSFINITE CONNECTEDNESS}
\author{A. H. Zemanian}
\date{}

\begin{document}
\maketitle
\baselineskip21pt

{\ Abstract --- A criterion is established for the transitivity
of connectedness in a transfinite graph.  Its proof is much 
shorter than a prior argument published previously for 
that criterion.
\\

Key Words: Transfinite graphs, transfinite connectedness, 
nondisconnectable tips.
} 

\vspace{.3in}

The objective of this work is provide a proof of 
a criterion for the transitivity of transfinite connectedness
that is much shorter than a previously published proof 
\cite[pages 57-63]{tgen} for that criterion.

We shall use the terminology defined in that book \cite{tgen}
and will add a few more definitions here. ${\cal G}^{\nu}$ will
denote a transfinite graph of rank $\nu$ with $0\leq\nu\leq\omega$.
Two branches in ${\cal G}^{\nu}$ are said to be $\rho$-{\em connected}
if there exists a $\gamma$-path $(0\leq\gamma\leq\rho\leq \nu)$ 
that terminates at a 0-node of each branch.  Then,
${\cal G}^{\nu}$ is called {\em transfinitely connected} if every two
of its branches are $\rho$-connected for some 
$\rho$ depending on the branches.  Transfinite connectedness 
is reflexive and symmetric binary relationship between the branches of 
${\cal G}^{\nu}$, but it need not be transitive.  Examples of such 
nontransitivity are given in \cite[Examples 3.1-5 and 3.1-6]{tgen}
and \cite[Example 3.1-1]{pg}.

However, the transitivity of $\rho$-connectedness
can be assured by restricting the 
graph ${\cal G}^{\nu}$ appropriately.  Two criteria, 
each of which suffices for such transitivity, are 
\cite[Condition 3.2-1 and Condition 3.5-1]{tgen}.  
The first criterion is rather 
lengthy and of limited use.  The second one is quite simple to state
and not overly restrictive,
but the proof of it given in \cite{tgen} is quite long and 
complicated.\footnote{Please
note the correction for \cite[pages 63-69]{tgen}
stated in the Errata for that book given in the URL:
www.ee.sunysb.edu/\~\,$\!$zeman. Namely, every elementary set of a sequence
of any rank is required to be node-distinct. This was implied by
the fact that the argument was based on paths,
but it should have been explicitly stated.} 
This latter case can be simplified by assuming that 
every node of ${\cal G}^{\nu}$ 
is {\em pristine}, that is, each node does not embrace a node of 
lower rank---and consequently
is not embraced by a node of higher 
rank either.\footnote{See \cite[Sec. 3.1]{pg}
in this regard.}   
Under this condition, the second criterion \cite[Condition 3.5-1]{tgen}
works again and with a much shorter proof.
It turns out that that latter proof can be
modified to apply to nonpristine graphs as well, yielding 
thereby a more elegant argument for the transitivity of
$\rho$-connectedness.  This is the subject of the present work.

We need some definitions. Two nonelementary 
tips (not necessarily of the the same rank) are said to be 
{\em disconnectable} if one can find two representatives, one 
for each tip, that are totally disjoint. (Any two tips, at least 
one of which is elementary, are 
simply taken to be {\em disconnectable}.)
As the negation of ``disconnectable,'' we say that two tips 
are {\em nondisconnectable} if every representative of one of them
meets every representative of the other tip infinitely often,
that is, at infinitely many nodes.
We can state this somewhat differently, as follows. 
Two tips, $t_{1}^{\gamma}$ and $t_{2}^{\delta}$, are called nondisconnectable
if $P_{1}^{\gamma}$ and $P_{2}^{\delta}$ meet at least once whenever
$P_{1}^{\gamma}$ is a representative of $t_{1}^{\gamma}$ and $P_{2}^{\delta}$
is a representative of $t_{2}^{\delta}$. 

Here is the sufficient condition that will ensure the transitivity
of connectedness.

{\bf Condition 1.}  {\em If two tips of ranks less than $\nu$
(and not necessarily of the 
same rank) are nondisconnectable, then either they are
shorted together or at least one of them is open.}

Before turning to our principle result, we need to introduce some more 
specificity into our terminology.  The maximal nodes 
partition all the nodes of ${\cal G}^{\nu}$;  
indeed, two nodes are in the same 
set of the partition if they are embraced by the same maximal 
node.
Thus, corresponding to any path, there is a unique set of maximal nodes 
that together embrace all the nodes of the path, and that set of maximal nodes 
becomes totally ordered in accordance with an orientation assigned to the
path.
Moreover, two maximal nodes are perforce totally disjoint.
In order to avoid explicating repeatedly which nodes embrace which nodes,
it is convenient to deal only with the maximal nodes
of our graph ${\cal G}^{\nu}$.

Let us expand upon this point because it is the key idea that allows
us to replace the long proof in \cite[pages 57-73]{tgen}
by a shorter, modified version of the proof in \cite[pages 31-35]{pg}.
Given any path $P^{\rho}$ of any rank $\rho$ 
in ${\cal G}^{\nu}$,
let ${\cal X}$ be the set of (possibly nonmaximal) nodes encountered
in the recursive construction of $P^{\rho}$;  that is, ${\cal X}$
consists of the nodes in the sequential representation of 
$P^{\rho}$ along with the nodes of lower ranks 
in the sequential representations of 
the paths between those nodes,
and also along with the nodes of still lower ranks 
in the sequential representations
of the paths of lower ranks in those latter paths, and so on down to 
the 0-nodes of the 0-paths embraced by $P^{\rho}$.
Each node of ${\cal X}$ need not be maximal, but it has a unique 
maximal node embracing it.  Let ${\cal X}_{m}$ denote the set of those
maximal nodes.  So, given $P^{\rho}$ and thereby ${\cal X}$ 
and ${\cal X}_{m}$, each node of ${\cal X}_{m}$ corresponds 
to a unique set of nodes 
in ${\cal X}$ such that, 
for every two nodes in that latter set, one embraces the other.
The union of those sets is ${\cal X}$.  In this way,
there is a bijection between ${\cal X}_{m}$ and the collection of the
said sets, which in fact comprise a partition of ${\cal X}$.
Moreover, when $P^{\rho}$ has an orientation, that orientation induces 
a total ordering of ${\cal X}_{m}$.  Thus, when dealing with the nodes
of $P^{\rho}$, we can fix our attention on the nodes of ${\cal X}_{m}$.  
This we shall do.  Any node $x^{\beta}$ of ${\cal X}_{m}$ will be called 
a {\em maximal node for} $P^{\rho}$ to distinguish it from
the corresponding subset of nodes in ${\cal X}$.  In general, we may have 
$\beta> \rho$, but not necessarily always.  Also, we can transfer 
our terminology for the nodes of a path, such as 
``incident to a node,'' ``meets a node,'' or ``terminates at a node,''
to those maximal nodes.
For example, a path is {\em incident to} a node $x\in {\cal X}_{m}$
if one of the path's terminal tips is embraced by $x$; this will 
be so if it has a terminal node embraced by $x$.  In the latter case, 
we say that the path {\em meets} x and also {\em terminates at} $x$.
Thus, two paths {\em meet at} $x\in{\cal X}_{m}$ if $x$ embraces a node $y$ 
of one path and a node $z$ of the second path such that 
$y$ embraces or is embraced by $z$.

The proof of our main result (Theorem 3) 
requires another idea, namely, ``path cuts.''  Let $P^{\rho}$
be a $\rho$-path with an orientation.  Let ${\cal Y}$ be the set 
of all branches and all (not necessarily maximal)
nodes of all ranks in $P^{\rho}$.  The orientation
of $P^{\rho}$ totally orders ${\cal Y}$.  With $y_{1}$ and $y_{2}$
being two members of ${\cal Y}$, we say that $y_{1}$ is {\em before} $y_{2}$
and that $y_{2}$ is {\em after} $y_{1}$ if in a tracing of $P^{\rho}$
in the direction of its orientation $y_{1}$ is met before $y_{2}$
is met. 
A {\em path cut} $\{{\cal B}_{1},{\cal B}_{2}\}$ 
for $P^{\rho}$ is a partitioning
of the set of branches of $P^{\rho}$ into two nonempty subsets, 
${\cal B}_{1}$ and ${\cal B}_{2}$, such that every branch of ${\cal B}_{1}$
is before every branch of ${\cal B}_{2}$.  Another way of stating this
is as follows.  The partition $\{{\cal B}_{1},{\cal B}_{2}\}$
of the branch set of $P^{\rho}$ comprises a path cut for $P^{\rho}$
if and only if, for each branch $b\in {\cal B}_{1}$, every branch 
of $P^{\rho}$ before $b$ is also a member of ${\cal B}_{1}$.

Here is a result that is easily proven through induction on ranks.

{\bf Lemma 2.}  {\em For each path cut $\{{\cal B}_{1},{\cal B}_{2}\}$
for $P^{\rho}$, there is a unique 
maximal node $x^{\gamma}$ $(\gamma\leq \rho)$
for $P^{\rho}$ such that every branch $b_{1}\in{\cal B}_{1}$ is before
$x^{\gamma}$ and every branch $b_{2}\in{\cal B}_{2}$ 
is after $x^{\gamma}$.}

We will say that the path cut {\em occurs at} the maximal node
$x^{\gamma}$. 
As was stated before, 
a node is called {\em maximal} if it is not embraced by a 
node of higher rank.  On the other hand, a maximal node 
for $P^{\rho}$ may 
embrace many other nodes.
It follows that all the nodes of $P^{\rho}$
other than the nodes of $P^{\rho}$ 
embraced by $x^{\gamma}$ 
are also partitioned into two sets,
the nodes of one set being before $x^{\gamma}$
and the nodes of the other set being after $x^{\gamma}$.

Here now is our main result.

{\bf Theorem 3.}  {\em Let ${\cal G}^{\nu}$ $(0\leq\nu\leq\omega)$ 
be a $\nu$-graph for which Condition 1 is satisfied.  Let 
$x_{a}$, $x_{b}$, and $x_{c}$ be three different maximal 
nodes in ${\cal G}^{\nu}$ such that, if any one of them is a singleton, 
its sole tip is disconnectable from every tip in the other two nodes.
If $x_{a}$ and $x_{b}$ are $\rho$-connected and if $x_{b}$
and $x_{c}$ are $\rho$-connected $(0\leq\rho\leq\nu)$,
then $x_{a}$ and $x_{c}$ are $\rho$-connected.}

{\bf Proof.}  Note that, if any one of $x_{a}$, $x_{b}$, and $x_{c}$ is a 
nonsingleton, then each of its embraced tips
(perforce, nonopen) must be disconnectable 
from every tip embraced by the other two nodes;  indeed, otherwise, 
the tips of that node would be shorted 
to the tips of another one of those three nodes, according to
Condition 1 and our hypothesis, 
and thus $x_{a}$, $x_{b}$, and $x_{c}$ could not be 
three different maximal nodes.

Let $P_{ab}^{\alpha}$ $(\alpha\leq\rho)$ be a two-ended
$\alpha$-path that terminates at 
the maximal nodes $x_{a}$ and $x_{b}$ and is oriented 
from $x_{a}$ to $x_{b}$, and let $P_{bc}^{\beta}$ $(\beta\leq\rho)$
be a two-ended $\beta$-path that terminates at 
the maximal nodes $x_{b}$ and $x_{c}$ and is
oriented from $x_{b}$ to $x_{c}$.  Let $P_{ba}^{\alpha}$ be 
$P_{ab}^{\alpha}$ but with the reverse orientation.  $P_{ba}^{\alpha}$
cannot have infinitely many $\alpha$-nodes because it is two-ended.

Let $\{ x_{i}\}_{i\in I}$ be the set of maximal nodes at which 
$P_{ba}^{\alpha}$ and $P_{bc}^{\beta}$ meet, and let ${\cal X}_{1}$ be that 
set of nodes with the order induced by the orientation 
of $P_{ba}^{\alpha}$.  If ${\cal X}_{1}$ has a last node $x_{l}$, then a 
tracing along $P_{ab}^{\alpha}$ from $x_{a}$ to $x_{l}$ followed by a tracing 
along $P_{bc}^{\beta}$ from $x_{l}$ to $x_{c}$ yields a path of rank no larger
than $\rho$ that connects $x_{a}$ and $x_{c}$.  
Thus, $x_{a}$ and $x_{c}$ are $\rho$-connected in this case.
This will certainly be so when $\{ x_{i}\}_{i\in I}$ is a finite set.

So, assume ${\cal X}_{1}$ is an infinite, ordered set (ordered as stated). 
We shall show that ${\cal X}_{1}$ has a last node.
Let $Q_{1}$ be the path induced by those branches of 
$P_{ba}^{\alpha}$ that lie between nodes of ${\cal X}_{1}$
(i.e., as $P_{ba}^{\alpha}$ is traced from $x_{b}$ onward, such a branch is 
traced after some node of ${\cal X}_{1}$ and before another node of 
${\cal X}_{1}$, those nodes depending upon the choice of the branch.)   
Let ${\cal B}_{1}$ be the set of those branches.
We can take it that $P^{\alpha}_{ba}$ extends beyond the nodes of 
${\cal X}_{1}$, for otherwise ${\cal X}_{1}$ would have $x_{a}$
as its last  node, and $x_{a}$ and $x_{c}$ would be $\rho$-connected.
Therefore,
we also have a nonempty set ${\cal B}_{2}$ 
consisting of those branches in $P_{ba}^{\alpha}$ that are not
in ${\cal B}_{1}$.  $\{ {\cal B}_{1},{\cal B}_{2}\}$ is a path cut for 
$P_{ba}^{\alpha}$.  Consequently, by Lemma 2., there is a unique
maximal node $x_{1}^{\gamma_{1}}$
at which that path cut occurs.
Thus, $Q_{1}$ terminates at $x_{1}^{\gamma_{1}}$. 
Let $t_{1}^{\rho_{1}}$
be the $\rho_{1}$-tip through which $Q_{1}$ reaches $x_{1}^{\gamma_{1}}$.
Every representative of $t_{1}^{\rho_{1}}$ contains infinitely 
many nodes of ${\cal X}_{1}$ (otherwise, ${\cal X}_{1}$
would have a last node).

Now, consider $P_{bc}^{\beta}$.  We can take it that there is a 
maximal node $x_{d}$ in $P_{bc}^{\beta}$ different from
$x_{c}$ such that the subpath of $P_{bc}^{\beta}$
between $x_{d}$ and $x_{c}$ is totally disjoint from $P_{ba}^{\alpha}$,
for otherwise $x_{a}$ and $x_{c}$ would have to be the same node
according to Condition 1 and our hypothesis again.
We can partition the branches of $P_{bc}^{\beta}$ into two
sets, ${\cal B}_{3}$ and ${\cal B}_{4}$, as follows.
Each branch of the first set ${\cal B}_{3}$ is 
such that it lies before (according to the orientation of 
$P_{bc}^{\beta}$) at least one node 
in ${\cal X}_{1}$ of each representative of 
$t_{1}^{\rho_{1}}$, this being so for all such representatives.\footnote{Let
us note here a correction for \cite[page 34, line 19 up]{pg}.
Replace the sentence on that line by the following sentence: 
``Let ${\cal N}_{2}$ be the set of those nodes in $\{ n_{i}\}_{i\in I}$
that lie before (according to the orientation of
$P_{bc}^{\beta}$) at least one node of $\{ n_{i}\}_{i \in I}$
in each representative of $t_{1}^{\rho_{1}}$, 
this being so for all such representatives.''}
The second set ${\cal B}_{4}$ consists of all the branches 
of $P_{bc}^{\beta}$ that are not in ${\cal B}_{3}$.
No branch of ${\cal B}_{4}$ can precede a branch of ${\cal B}_{3}$.
Thus, we have a path cut $\{{\cal B}_{3},{\cal B}_{4}\}$
for $P_{bc}^{\beta}$ and thereby (according to Lemma 2.)
a unique maximal node $x_{2}^{\gamma_{2}}$ lying after the branches of 
${\cal B}_{3}$ and before the branches of ${\cal B}_{4}$.
Let $Q_{2}$ be the path induced by ${\cal B}_{3}$.  It reaches 
$x_{2}^{\gamma_{2}}$ through some tip $t_{2}^{\rho_{2}}$.
Furthermore, each representative of $t_{2}^{\rho_{2}}$ must meet 
each representative of $t_{1}^{\rho_{1}}$ at least once
because they meet at at least one node of 
${\cal X}_{1}$.  Thus, $t_{1}^{\rho_{1}}$ and $t_{2}^{\rho_{2}}$
are nondisconnectable. Moreover, neither of those tips can be open 
(i.e., be in a singleton node) because the paths $P_{ba}^{\alpha}$
and $P_{bc}^{\beta}$ pass through and beyond their respective nodes 
$x_{1}^{\gamma_{1}}$ and $x_{2}^{\gamma_{2}}$.  So, by Condition 1, 
$x_{1}^{\gamma_{1}}$ and $x_{2}^{\gamma_{2}}$ must be the same 
node because the tips $t_{1}^{\rho_{1}}$ and $t_{2}^{\rho_{2}}$
are shorted together.  

This means that ${\cal X}_{1}$ has a last node, namely,
$x_{1}^{\gamma_{1}}= x_{2}^{\gamma_{2}}$.  It follows now that 
$x_{a}$ and $x_{c}$ are $\rho$-connected. $\clubsuit$

The last proof has established the following two result.

{\bf Corollary 4.}  {\em Under the hypothesis of Theorem 3,
let $P_{ab}^{\alpha}$ be a two-ended
$\alpha$-path connecting nodes $x_{a}$ and $x_{b}$, and let 
$P_{bc}^{\beta}$ be a two-ended $\beta$-path connecting 
nodes $x_{b}$ and $x_{c}$.
Then, there is a two-ended $\gamma$-path ($\gamma\leq \max\{\alpha,\beta\}$)
connecting $x_{a}$ and $x_{c}$ that lies in 
$P_{ab}^{\alpha}\cup P_{bc}^{\beta}$.}

{\bf Corollary 5.}  {\em Under the hypothesis of Theorem 3, let 
$\{ x_{i}\}_{i\in I}$ be the set of maximal nodes at which 
two two-ended paths meet.  Assume that, if any terminal node
of either path is a singleton, its sole tip is disconnectable 
from every tip in the two terminal nodes of the other path.
Assign to $\{ x_{i}\}_{i\in I}$ the total ordering 
induced by an orientation of one of these two paths.  Then, 
$\{ x_{i}\}_{i\in I}$ has both a first node and a last node.}

\end{document}